\documentclass[nospthms]{svmult}

\usepackage[numbers]{natbib}

\newcommand{\co}{:}

\usepackage{verbatim,amssymb,amsmath,amsthm,url,cite}

\usepackage{bbm,stmaryrd}

\usepackage[cmtip,arrow]{xy}

\newtheorem{theorem}{Theorem}[section]
\newtheorem*{theorem*}{Theorem}
\newtheorem{proposition}[theorem]{Proposition}
\newtheorem*{proposition*}{Proposition}
\newtheorem{corollary}[theorem]{Corollary}
\newtheorem{lemma}[theorem]{Lemma}
\theoremstyle{definition}

\newtheorem{definition}[theorem]{Definition}
\newtheorem*{definition*}{Definition}
\newtheorem*{remark*}{Remark}

\input diagxy
\def\too{\ifnextchar/{\toop}{\toop/>/}}%

\renewcommand{\1}{\mathbbm{1}}

\newcommand{\adjunction}[4]{#1 : \, #2 \rightleftharpoons #3 \, : #4}

\renewcommand{\[}{\llbracket}
\renewcommand{\]}{\rrbracket}
\newcommand{\N}{\mathbb{N}}
\newcommand{\RR}{\mathbb{R}}

\newcommand{\p}{\partial}
\newcommand{\Z}{\mathbb{Z}}

\newcommand{\C}{\mathcal{C}}

\renewcommand{\E}{\mathcal{E}}
\newcommand{\V}{\mathcal{V}}
\renewcommand{\o}{\otimes}

\DeclareMathOperator{\Ob}{Ob}

\DeclareMathOperator{\sgn}{sgn}
\DeclareMathOperator{\Diff}{Diff}

\newcommand{\Cat}{\mathsf{Cat}}
\newcommand{\Set}{\mathsf{Set}}
\newcommand{\sSet}{\mathsf{sSet}}
\newcommand{\VCat}{\text{$\V$-$\mathsf{Cat}$}}

\renewcommand{\P}{\mathcal{P}}
\newcommand{\Q}{\mathcal{Q}}
\newcommand{\R}{\mathcal{R}}

\renewcommand{\)}{)\!)}
\newcommand{\<}{\langle}
\renewcommand{\>}{\rangle}

\newcommand{\Cob}{\mathsf{Cob}}
\newcommand{\cob}{\mathsf{cob}}

\newcommand{\CA}{\mathcal{A}}
\newcommand{\CB}{\mathcal{B}}
\newcommand{\CN}{\mathcal{N}}
\newcommand{\CP}{\mathcal{P}}
\newcommand{\CS}{\mathsf{S}}
\newcommand{\CT}{\mathsf{T}}
\newcommand{\ct}{\mathsf{t}}
\newcommand{\K}{\mathsf{s}}
\newcommand{\TT}{\mathbb{T}}

\renewcommand{\SS}{\mathbb{S}}

\newcommand{\op}{\circ}
\newcommand{\tint}{{\textstyle\int}}
\newcommand{\presheaf}[1]{#1{}\hat{\ }}

\newcommand{\h}{\mathsf{h}}
\DeclareMathOperator{\Mod}{Mod}
\DeclareMathOperator{\PreMod}{PreMod}

\newcommand{\SSS}{\SS\wr}

\newcommand{\CC}{\mathbb{C}}

\DeclareMathOperator{\PreOp}{PreOp}
\DeclareMathOperator{\Op}{Op}
\DeclareMathOperator{\End}{\mathsf{End}}
\DeclareMathOperator{\Hom}{Hom}
\DeclareMathOperator{\tr}{tr}

\DeclareMathOperator{\SO}{SO}

\def\to{\mathchoice{\longrightarrow}{\rightarrow}{\rightarrow}{\rightarrow}}
\def\To{\mathchoice{\Longrightarrow}{\Rightarrow}{\Rightarrow}{\Rightarrow}}

\makeatletter
\def\rightharpoonfill@{\arrowfill@\relbar\relbar\rightharpoonup}
\newcommand{\overrightharpoon}{\mathpalette{\overarrow@\rightharpoonfill@}}
\def\hookrightarrow{\mathchoice
  {\DOTSB\lhook\joinrel\relbar\joinrel\rightarrow}
  {\DOTSB\lhook\joinrel\rightarrow}
  {\DOTSB\lhook\joinrel\rightarrow}
  {\DOTSB\lhook\joinrel\rightarrow}}
\makeatother

\newcommand{\CG}{\mathcal{G}}
\newcommand{\CH}{\mathcal{H}}
\newcommand{\CD}{\mathcal{D}}
\DeclareMathOperator{\vCob}{\Cob^\rightharpoonup}
\DeclareMathOperator{\vcob}{\cob^\rightharpoonup}

\renewcommand{\phi}{\varphi}
\renewcommand{\*}{\cdot}

\newcommand{\reflex}{\three/>`<-`>/}

\begin{document}

\title*{Operads revisited}

\author{Ezra Getzler}

\institute{Northwestern University, Evanston, Illinois, USA \\
\texttt{getzler@northwestern.edu}}

\maketitle

% \primaryclass{18D50, 57R56}

% \secondaryclass{18D10, 18D20, 18C35}

%\keywords{Operads, symmetric monoidal categories, enriched
% categories, topological field theory, locally presentable categories}

% \cl{\emph{To Yuri Manin, many happy returns.}}

This paper presents an approach to operads related to recent work in
topological field theory (\citet{Costello}).  The idea is to represent
operads as symmetric monoidal functors on a symmetric monoidal
category $\CT$; we recall how this works for cyclic and modular
operads and dioperads in Section~\ref{discrete}. This point of view
permits the construction of all sorts of variants on the notion of an
operad. As our main example, we present a simplicial variant of
modular operads, related to Segal's definition of quantum field
theory, as modified for topological field theory (\citet{bv}). (This
definition is only correct in a cocomplete symmetric monoidal category
whose tensor product preserves colimits, but this covers most cases of
interest.)

An operad in the category $\Set$ of sets may be presented as a
symmetric monoidal category $\ct\P$, called the theory associated to
$\P$ (\citet{BoardmanVogt}).  The category $\ct\P$ has the natural
numbers as its objects; tensor product is given by addition. The
morphisms of $\ct\P$ are built using the operad $\P$:
\begin{equation*}
  \ct\P(m,n) = \bigsqcup_{f\co\{1,\dotsc,m\}\to\{1,\dotsc,n\}}
  \bigsqcap_{i=1}^n \P(|f^{-1}(i)|) .
\end{equation*}
The category of $\P$-algebras is equivalent to the category of
symmetric monoidal functors from $\ct\P$ to $\Set$; this reduces the
study of algebras over operads to the study of symmetric monoidal
functors.

The category of contravariant functors from the opposite category
$\CT^\op$ of a small category $\CT$ to the category $\Set$ of sets
\begin{equation*}
  \presheaf{\CT} = [\CT^\op,\Set]
\end{equation*}
is called the category of presheaves of $\CT$. If $\CT$ is a symmetric
monoidal category, then $\presheaf{\CT}$ is too (\citet{Day}); its
tensor product is the coend
\begin{equation*}
  V \ast W = \tint^{A,B\in\CT} \CT(-,A\o B) \times V(A) \times W(B) .
\end{equation*}

A symmetric monoidal functor $F \co \CS \to \CT$ between symmetric
monoidal categories $\CS$ and $\CT$ is a functor $F$ together with a
natural equivalence
\begin{equation*}
  \Phi \co \o \circ F\times F \To F \circ \o .  
\end{equation*}
The functor $F$ is lax symmetric monoidal if $\Phi$ is only a natural
transformation.

If $\tau\co\CS\to\CT$ is a symmetric monoidal functor, it is not
always the case that the induced functor
\begin{equation*}
  \presheaf{\tau} \co \presheaf{\CT} \to \presheaf{\CS}
\end{equation*}
is a symmetric monoidal functor; in general, it is only a lax
symmetric monoidal functor, with respect to the natural transformation
\begin{equation*}
  \begin{xy}
    \morphism(0,0)/=/<1500,0>[\presheaf{\tau}(V) \ast
    \presheaf{\tau}(W)`\tint^{A,B\in\CS} \CS(-,A\o B) \times V(\tau A)
    \times W(\tau B);]
    \morphism(0,0)/=>/<0,-1500>[\presheaf{\tau}(V)
    \ast \presheaf{\tau}(W)`\presheaf{\tau}(V\ast W);\Phi_{V,W}]
    \morphism(1500,0)|r|/=>/<0,-500>[\tint^{A,B\in\CS} \CS(-,A\o B)
    \times V(\tau A)\times W(\tau B)`\tint^{A,B\in\CS} \CT( \tau(-) ,
    \tau(A\o B))\times V(\tau A) \times W(\tau B);
    \text{ $\tau$ is a functor}]
    \morphism(1500,-500)|r|/=>/<0,-500>[\tint^{A,B\in\CS} \CT( \tau(-)
    , \tau(A\o B))\times V(\tau A) \times W(\tau B)`\tint^{A,B\in\CS}
    \CT( \tau(-) , \tau(A)\o\tau(B)) \times V(\tau A) \times W(\tau
    B); \text{ $\tau$ is symmetric monoidal}]
    \morphism(1500,-1000)|r|/=>/<0,-500>[\tint^{A,B\in\CS} \CT(
    \tau(-) , \tau(A)\o\tau(B) ) \times V(\tau A) \times W(\tau B)`
    \tint^{A,B\in\CT} \CT( \tau(-) , A\o B ) \times V(A) \times W(B) ;
    \text{ universality of coends}]
    \morphism(0,-1500)/=/<1500,0>[\presheaf{\tau}(V\ast W) `
    \tint^{A,B\in\CT} \CT( \tau(-) , A\o B ) \times V(A) \times W(B) ;
    ]
  \end{xy}
\end{equation*}

The following definition introduces the main object of study of this
paper.
\begin{definition*}
  A \textup{\textbf{pattern}} is a symmetric monoidal functor
  $\tau:\CS\to\CT$ between small symmetric monoidal categories $\CS$
  and $\CT$ such that
  $\presheaf{\tau}:\presheaf{\CT}\to\presheaf{\CS}$ is a symmetric
  monoidal functor (in other words, the natural transformation $\Phi$
  defined above is an equivalence).
\end{definition*}

Let $\tau$ be a pattern and let $\C$ be a symmetric monoidal category.
A \textbf{$\tau$-preoperad} in $\C$ is a symmetric monoidal functor
from $\CS$ to $\C$, and a \textbf{$\tau$-operad} in $\C$ is a
symmetric monoidal functor from $\CT$ to $\C$. Denote by
$\PreOp_\tau(\C)$ and $\Op_\tau(\C)$ the categories of
$\tau$-preoperads and $\tau$-operads respectively. If $\C$ is a
cocomplete symmetric monoidal category, there is an adjunction
\begin{equation*}
  \adjunction{\tau_*}{\PreOp_\tau(\C)}{\Op_\tau(\C)}{\tau^*} .
\end{equation*}
In Section 2, we prove that if $\tau$ is essentially surjective, $\C$
is cocomplete, and the functor $A\o B$ on $\C$ preserves colimits in
each variable, then the functor $\tau^*$ is monadic. We also prove
that if $\CS$ is a free symmetric monoidal category and $\C$ is
locally finitely presentable, then $\Op_\tau(\C)$ is locally finitely
presentable.

Patterns generalize the coloured operads of \citet{BoardmanVogt},
which are the special case where $\CS$ is the free symmetric category
generated by a discrete category (called the set of colours). Operads
are themselves algebras for a coloured operad, whose colours are the
natural numbers (cf.\ \citet{BM2}), but it is more natural to think of
the colour $n$ as having nontrivial automorphisms, namely the
symmetric group $\SS_n$.

The definition of a pattern may be applied in the setting of
simplicial categories, or even more generally, enriched categories. We
define and study patterns enriched over a symmetric monoidal category
$\V$ in Section~\ref{pattern}; we recall those parts of the formalism
of enriched categories which we will need in the appendix.

In Section 3, we present examples of a simplicial pattern which arises
in topological field theory. The most interesting of these is related
to modular operads. Let $\cob$ be the (simplicial) groupoid of
diffeomorphisms between compact connected oriented surfaces with
boundary. Let $S_{g,n}$ be a compact oriented surface of genus $g$
with $n$ boundary circles, and let the mapping class group be the
group of components of the oriented diffeomorphism group:
\begin{equation*}
  \Gamma_{g,n} = \pi_0( \Diff_+(S_{g,n}) ) .
\end{equation*}
The simplicial groupoid $\cob$ has a skeleton
\begin{equation*}
  \bigsqcup_{g,n} \, \Diff_+(S_{g,n}) ,
\end{equation*}
and if $2g-2+n>0$, there is a homotopy equivalence $\Diff_+(S_{g,n}) =
\Gamma_{g,n}$. In Section~\ref{cob}, we define a pattern whose
underlying functor is the inclusion $\SSS\cob\hookrightarrow\Cob$.
This pattern bears a similar relation to modular operads that braided
operads (\citet{braided}) bear to symmetric operads. The mapping class
group $\Gamma_{0,n}$ is closely related to the ribbon braid group
$B_n\wr\Z$. Denoting the generators of $B_n\subset B_n\wr\Z$ by
$\{b_1,\dotsc,b_{n-1}\}$ and the generators of $\Z^n\subset B_n\wr\Z$
by $\{t_1,\dotsc,t_n\}$, Moore and Seiberg show \citep[Appendix
B.1]{MS} that $\Gamma_{0,n}$ is isomorphic to the quotient of
$B_n\wr\Z$ by its subgroup $\< (b_1b_2\dotsc b_{n-1}t_n)^n ,
b_1b_2\dots b_{n-1}t_n^2b_{n-1}\dots b_1 \>$. In this sense, operads
for the pattern $\Cob$ are a cyclic analogue of braided operads.

\subsection*{Acknowledgements}

I thank Michael Batanin and Steve Lack for very helpful feedback on
earlier drafts of this paper.

I am grateful to a number of institutions for hosting me during the
writing of this paper: RIMS, University of Nice, and KITP.  I am
grateful to my hosts at all of these institutions, Kyoji Saito and
Kentaro Hori at RIMS, Andr\'e Hirschowitz and Carlos
Simpson in Nice, and David Gross at KITP, for their hospitality.

I am also grateful to Jean-Louis Loday for the opportunity to lecture
on an earlier version of this work at Luminy. I received support from
NSF Grants DMS-0072508 and DMS-0505669.

\section{Modular operads as symmetric monoidal functors} \label{discrete}

In this section, we define modular operads in terms of the symmetric
monoidal category of dual graphs. Although we do not assume
familiarity with the original definition (\citet{modular}; see also
\citet{MSS}), this section will certainly be easier to understand if
this is not the reader's first brush with the subject.

We also show how modifications of this construction, in which dual
graphs are replaces by forests, or by directed graphs, yielding cyclic
operads and modular dioperads. Finally, we review the definition of
algebras for modular operads and dioperads.

\subsection*{Graphs}

A graph $\Gamma$ consists of the following data:
\begin{enumerate}
\item finite sets $V(\Gamma)$ and $F(\Gamma)$, the sets of vertices
  and flags of the graph;
\item a function $p\co F(\Gamma)\to V(\Gamma)$, whose fibre
  $p^{-1}(v)$ is the set of flags of the graph meeting at the vertex
  $v$;
\item an involution $\sigma\co F(\Gamma)\to F(\Gamma)$, whose
  fixed points are called the legs of $\Gamma$, and whose remaining
  orbits are called the edges of $\Gamma$.
\end{enumerate}
We denote by $L(\Gamma)$ and $E(\Gamma)$ the sets of legs and edges of
$\Gamma$, and by $n(v)=|p^{-1}(v)|$ the number of flags meeting a
vertex $v$.

To a graph is associated a one-dimensional cell complex, with 0-cells
$V(\Gamma)\sqcup L(\Gamma)$, and 1-cells $E(\Gamma)\sqcup L(\Gamma)$.
The 1-cell associated to an edge $e=\{f,\sigma(f)\}$, $f\in F(G)$, is
attached to the 0-cells corresponding to the vertices $p(f)$ and
$p(\sigma(f))$ (which may be equal), and the 1-cell associated to a
leg $f\in L(G)\subset F(G)$ is attached to the 0-cells corresponding
to the vertex $p(f)$ and the leg $f$ itself.

The edges of a graph $\Gamma$ define an equivalence relation on its
vertices $V(\Gamma)$; the components of the graph are the equivalence
classes with respect to this relation. Denote the set of components by
$\pi_0(\Gamma)$. The Euler characteristic of a component $C$ of a
graph is $e(C)=|V(C)|-|E(C)|$. Denote by $n(C)$ the number of legs of
a component $C$.

\subsection*{Dual graphs}
A dual graph is a graph $\Gamma$ together with a function
$g\co V(\Gamma)\to\N$. The natural number $g(v)$ is called the
genus of the vertex $v$.

The genus $g(C)$ of a component $C$ of a dual graph $\Gamma$ is
defined by the formula
\begin{equation*}
  g(C) = \sum_{v\in C} g(v) + 1 - e(C) .
\end{equation*}
The genus of a component is a non-negative integer, which equals $0$
if and only if $C$ is a tree and $g(v)=0$ for all vertices $v$ of $C$.

A \textbf{stable} graph is a dual graph $\Gamma$ such that for all
vertices $v\in V(\Gamma)$, the integer $2g(v)-2+n(v)$ is positive. In
other words, if $g(v)=0$, then $n(v)$ is at least $3$, while if
$g(v)=1$, then $n(v)$ is nonzero. If $\Gamma$ is a stable graph, then
$2g(C)-2+|L(C)|>0$ for all components $C$ of $\Gamma$.

\subsection*{The symmetric monoidal category $\CG$}

The objects of the category $\CG$ are the dual graphs $\Gamma$ whose
set of edges $E(\Gamma)$ is empty.  Equivalently, an object of $\CG$
is a pair of finite sets $L$ and $V$ and functions $p\co L\to V$,
$g\co V\to\N$.

A morphism of $\CG$ with source $(L_1,V_1,p_1,g_1)$ and target
$(L_2,V_2,p_2,g_2)$ is a pair consisting of a dual graph $\Gamma$ with
$F(\Gamma)=L_1$, $V(\Gamma)=V_1$, and $g(\Gamma)=g_1$, together with
isomorphisms $\alpha\co L_2\to L(\Gamma)$ and
$\beta\co V_2\to\pi_0(\Gamma)$ such that $p\circ\alpha=\beta\circ
p_2$ and $\alpha^*g=g_2$.

The definition of the composition of two morphisms
$\Gamma=\Gamma_2\circ\Gamma_1$ in $\CG$ is straightforward. We have
dual graphs $\Gamma_1$ and $\Gamma_2$, with
\begin{align*}
  F(\Gamma_1) &= L_1 , & V(\Gamma_1) &= V_1 , & g(\Gamma_1) &= g_1 , \\
  F(\Gamma_2) &= L_2 , & V(\Gamma_2) &= V_2 , & g(\Gamma_2) &= g_2 ,
\end{align*}
together with isomorphisms
\begin{align*}
  \alpha_1\co L_2 &\to L(\Gamma_1) , & \beta_1\co V_2 &\to \pi_0(\Gamma_1) , \\
  \alpha_2\co L_3 &\to L(\Gamma_2) , & \beta_2\co V_3 &\to \pi_0(\Gamma_2) .
\end{align*}

Since the source of $\Gamma=\Gamma_2\circ\Gamma_1$ equals the source
of $\Gamma_1$, we see that $p\co F(\Gamma)\to V(\Gamma)$ is
identified with $p\co L_1\to V_1$. The involution $\sigma$ of
$L_1$ is defined as follows: if $f\in F(\Gamma)=L_1$ lies in an edge
of $\Gamma_1$, then $\sigma(f)=\sigma_1(f)$, while if $f$ is a leg of
$\Gamma_1$, then $\sigma(f)=\alpha_1(\sigma_2(\alpha_1^{-1}(f)))$. The
isomorphisms $\alpha$ and $\beta$ of $\Gamma$ are simply the
isomorphisms $\alpha_2$ and $\beta_2$ of $\Gamma_2$.

In other words, $\Gamma_2\circ\Gamma_1$ is obtained from $\Gamma_1$ by
gluing those pairs of legs of $\Gamma_1$ together which correspond to
edges of $\Gamma_2$. It is clear that composition in $\CG$ is
associative.

The tensor product on objects of $\CG$ extends to morphisms, making
$\CG$ into a symmetric monoidal category.

A morphism in $\CG$ is invertible if the underlying stable graph has
no edges, in other words, if it is simply an isomorphism between two
objects of $\CG$. Denote by $\CH$ the groupoid consisting of all
invertible morphisms of $\CG$, and by $\tau:\CH\hookrightarrow\CG$ the
inclusion. The groupoid $\CH$ is the free symmetric monoidal functor
generated by the groupoid $\h$ consisting of all morphisms of $\CG$
with connected domain.

The category $\CG$ has a small skeleton; for example, the full
subcategory of $\CG$ in which the sets of flags and vertices of the
objects are subsets of the set of natural numbers. The tensor product
of $\CG$ takes us outside this category, but it is not hard to define
an equivalent tensor product for this small skeleton. We will tacitly
replace $\CG$ by this skeleton, since some of our constructions will
require that $\CG$ be small.

\subsection*{Modular operads as symmetric monoidal functors}

The following definition of modular operads may be found in
\citet{Costello}. Let $\C$ be a symmetric monoidal category which is
cocomplete, and such that the functor $A\o B$ preserves colimits in
each variable. (This last condition is automatic if $\C$ is a closed
symmetric monoidal category.) A \textbf{modular preoperad} in $\C$ is
a symmetric monoidal functor from $\CH$ to $\C$, and a \textbf{modular
  operad} in $\C$ is a symmetric monoidal functor from $\CG$ to $\C$.

Let $\CA$ be a small symmetric monoidal category, and let $\CB$ be
symmetric monoidal category. Denote by $[\CA,\CB]$ the category of
functors and natural equivalences from $\CA$ to $\CB$, and by
$\[\CA,\CB\]$ the category of symmetric monoidal functors and monoidal
natural equivalences.

With this notation, the categories of modular operads, respectively
preoperads, in a symmetric monoidal category $\C$ are
$\Mod(\C)=\[\CG,\C\]$ and $\PreMod(\C)=\[\CH,\C\]$.
\begin{theorem*}
  \begin{enumerate}
  \item There is an adjunction
    \begin{equation*}
      \adjunction{\tau_*}{\PreMod(\C)}{\Mod(\C)}{\tau^*} ,
    \end{equation*}
    where $\tau^*$ is restriction along
    $\tau\co\CH\hookrightarrow\CG$, and $\tau_*$ is the coend
    \begin{equation*}
      \tau_*\P = \tint^{A\in\CH} \CG(A,-) \times \P(A) .
    \end{equation*}
  \item The functor
    \begin{equation*}
      \tau^* \co \Mod(\C) \to \PreMod(\C)
    \end{equation*}
    is monadic. That is, there is an equivalence of categories
    \begin{equation*}
      \Mod(\C) \simeq \PreMod(\C)^\TT ,
    \end{equation*}
    where $\TT$ is the monad $\tau^*\tau_*$.    
  \item The category $\Mod(\C)$ is locally finitely presentable if
    $\C$ is locally finitely presentable.
  \end{enumerate}
\end{theorem*}

In the original definition of modular operads (\citet{modular}), there
was an additional stability condition, which may be phrased in the
following terms. Denote by $\CG_+$ the subcategory of $\CG$ consisting
of stable graphs, and let $\CH_+=\CH\cap\CG_+$. Then a stable modular
preoperad in $\C$ is a symmetric monoidal functor from $\CH_+$ to
$\C$, and a stable modular operad in $\C$ is a symmetric monoidal
functor from $\CG_+$ to $\C$. These categories of stable modular
preoperads and operads are equivalent to the categories of stable
$\SS$-modules and modular operads of loc.\ cit.

\subsection*{Algebras for modular operads}

To any object $M$ of a closed symmetric monoidal category $\C$ is
associated the monoid $\textup{End}(M)=[M,M]$: an $A$-module is a
morphism of monoids $\rho\co A\to\textup{End}(M)$. The analogue of this
construction for modular operads is called a $\P$-algebra.

A bilinear form with domain $M$ in a symmetric monoidal category $\C$
is a morphism $t\co M\o M \to \1$ such that $t\circ\sigma=t$. Associated
to a bilinear form $(M,t)$ is a modular operad $\End(M,t)$, defined on
a connected object $\Gamma$ of $\CG$ to be
\begin{equation*}
  \End(M,t)(\Gamma) = M^{\o L(\Gamma)} .
\end{equation*}
If $\P$ is a modular operad, a $\P$-algebra $M$ is an object $M$ of
$\C$, a bilinear form $t\co M\o M\to \1$, and a morphism of modular
operads $\rho\co\P\to\End(M,t)$.

This modular operad has an underlying stable modular operad
$\End_+(M,t)$, defined by restriction to the stable graphs in $\CG$.

\subsection*{Cyclic operads}
A variant of the above definition of modular operads is obtained by
taking the subcategory of forests $\CG_0$ in the category $\CG$: a
forest is a graph each component of which is simply connected.
Symmetric monoidal functors on $\CG_0$ are cyclic operads.

Denote the cyclic operad underlying $\End(M,t)$ by $\End_0(M,t)$: thus
\begin{equation*}
  \End_0(M,t)(\Gamma) = M^{\o L(\Gamma)} .
\end{equation*}
If $\P$ is a cyclic operad, a $\P$-algebra $M$ is bilinear form $t$
with domain $M$ and a morphism of cyclic operads
$\rho\co\P\to\End_0(M,t)$.

There is a functor $\P\mapsto\P_0$, which associates to a modular
operad its underlying cyclic operad: this is the restriction functor
from $\[\CG,\C\]$ to $\[\CG_0,\C\]$.

There is also a stable variant of cyclic operads, in which $\CG_0$ is
replaced by its stable subcategory $\CG_{0+}$, defined by restricting
to forests in which each vertex meets at least three flags.

\subsection*{Dioperads and modular dioperads}

Another variant of the definition of modular operads is obtained by
replacing the graphs $\Gamma$ in the definition of modular operads by
digraphs (directed graphs):

A \textbf{digraph} $\Gamma$ is a graph together with a partition
\begin{equation*}
  F(\Gamma) = F_+(\Gamma) \sqcup F_-(\Gamma)
\end{equation*}
of the flags into outgoing and incoming flags, such that each edge has
one outgoing and one incoming flag. Each edge of a digraph has an
orientation, running towards the outgoing flag. The set of legs of a
digraph are partitioned into the outgoing and incoming legs:
$L_\pm(\Gamma)=L(\Gamma)\cap F_\pm(\Gamma)$. Denote the number of
outgoing and incoming legs by $n_\pm(\Gamma)$.

A dual digraph is a digraph together with a function $g\co
V(\Gamma)\to\N$. Imitating the construction of the symmetric monoidal
category of dual graphs $\CG$, we may construct a symmetric monoidal
category of dual digraphs $\CD$. A \textbf{modular dioperad} is a
symmetric monoidal functor on $\CD$. (These are the wheeled props
studied in a recent preprint of \citet{Merkulov}.)

If $M_\pm$ are objects of $\C$ and $t\co M_+\o M_-\to\1$ is a pairing, we
may construct a modular dioperad $\End^\rightharpoonup(M_\pm,t)$,
defined on a connected object $\Gamma$ of $\CD$ to be
\begin{equation}
  \label{modulardioperad}
  \End^\rightharpoonup(M)(\Gamma) = M_+^{\o  L_+(\Gamma)} \o
  M_-^{\o L_-(\Gamma)} .
\end{equation}
If $\P$ is a modular dioperad, a $\P$-algebra is a pairing $t\co M_+\o
M_-\to\1$ and a morphism of modular dioperads
$\rho\co\P\to\End^\rightharpoonup(M_\pm,t)$.

A directed forest is a directed graph each component of which is
simply connected; let $\CD_0$ be the subcategory of $\CD$ of
consisting of directed forests. A \textbf{dioperad} is a symmetric
monoidal functor on $\CD_0$ (\citet{dioperad}).

If $M$ is an object of $\C$, we may construct a dioperad
$\End^\rightharpoonup_0(M)$, defined on a connected object $\Gamma$ of
$\CD_0$ to be
\begin{equation*}
  \End^\rightharpoonup_0(M)(\Gamma) = \Hom\bigl( M^{\o
    L_+(\Gamma)},M^{\o L_-(\Gamma)} \bigr) .
\end{equation*}
When $M=M_-$ is a rigid object with dual $M^\vee=M_+$, this is a
special case of \eqref{modulardioperad}. If $\P$ is a dioperad, a
$\P$-algebra $M$ is an object $M$ of $\C$ and a morphism of dioperads
$\rho\co\P\to\End^\rightharpoonup_0(M)$.

\subsection*{props}
MacLane's notion \citep{MacLane} of a prop also fits into the above
framework. There is a subcategory $\CD_{\mathsf{P}}$ of $\CD$,
consisting of all dual digraphs $\Gamma$ such that each vertex has
genus $0$, and $\Gamma$ has no directed circuits. A prop is a
symmetric monoidal functor on $\CD_{\mathsf{P}}$. (This follows from
the description of free props in \citet{EE}.) Note that $\CD_0$ is a
subcategory of $\CD_{\mathsf{P}}$: thus, every prop has an underlying
dioperad. Note also that the dioperad $\End^\rightharpoonup_0(M)$ is
in fact a prop; if $\P$ is a prop, we may define a $\P$-algebra $M$ to
be an object $M$ of $\C$ and a morphism of props
$\rho\co\P\to\End^\rightharpoonup_0(M)$.

\section{Patterns} \label{pattern}

Patterns abstract the approach to modular operads sketched in
Section~1. This section develops the theory of patterns enriched over
a complete, cocomplete, closed symmetric monoidal category $\V$. In
fact, we are mainly interested in the cases where $\V$ is the category
$\Set$ of sets or the category $\sSet$ of simplicial sets. We refer to
the appendices for a review of the needed enriched category theory.

\subsection*{Symmetric monoidal $\V$-categories}

Let $\CA$ be a $\V$-category. Denote by $\SS_n$ the symmetric group on
$n$ letters. The wreath product $\SS_n\wr\CA$ is the $\V$-category
\begin{equation*}
  \SS_n\wr\CA = \SS_n \times \CA^n .
\end{equation*}
If $\alpha,\beta\in\SS_n$, the composition of morphisms
$(\alpha,\phi_1,\dotsc,\phi_n)$ and $(\beta,\psi_1,\dotsc,\psi_n)$ is
\begin{equation*}
  (\beta\circ\alpha,\psi_{\alpha_1}\circ\phi_1,\dotsc,
  \psi_{\alpha_n}\circ\phi_n) .  
\end{equation*}
Define the wreath product $\SSS\CA$ to be
\begin{equation*}
  \SSS\CA = \bigsqcup_{n=0}^\infty \SS_n\wr\CA .
\end{equation*}

In fact, $\SSS{(-)}$ is a 2-functor from the 2-category $\VCat$ to
itself. This 2-functor underlies a 2-monad $\SSS{(-)}$ on $\VCat$: the
composition
\begin{equation*}
  m\co\SSS{\SSS{(-)}}\to\SSS{(-)}
\end{equation*}
is induced by the natural inclusions
\begin{equation*}
  (\SS_{n_1}\times\CA^n_1)\times\dotsm\times(\SS_{n_k}\times\CA^n_k)
  \hookrightarrow \SS_{n_1+\dotsb+n_k} \times \CA^{n_1+\dotsb+n_k} ,
\end{equation*}
and the unit $\eta\co1_\VCat\to\SSS{(-)}$ is induced by the
natural inclusion
\begin{equation*}
  \CA \cong \SS_1\times\CA\hookrightarrow\SSS\CA .
\end{equation*}

\begin{definition}
  A \textbf{symmetric monoidal} $\V$-category $\C$ is a pseudo
  $\SSS{(-)}$-algebra in $\VCat$.
\end{definition}

If $\C$ is a symmetric monoidal category, we denote the object
obtained by acting on the object $(A_1,\dotsc,A_n)$ of $\SS_n\wr\C$ by
$A_1\o\dotsb\o A_n$.  When $n=0$, we obtain an object $\1$ of $\C$,
called the identity.  When $n=1$, we obtain a $\V$-endofunctor of
$\C$, which is equivalent by the natural $\V$-equivalence $\iota$ in
the definition of a pseudo $\SS_n\wr\C$-algebra to the identity
$\V$-functor. When $n=2$, we obtain a $\V$-functor $(A,B)\mapsto A\o
B$ from $\C\times\C$ to $\C$, called the tensor product. Up to
$\V$-equivalence, all of the higher tensor products are obtained by
iterating the tensor product $A\o B$: if $n>2$, there is a natural
$\V$-equivalence between the functors $A_1\o\dotsb\o A_n$ and
$(A_1\o\dotsb\o A_{n-1})\o A_n$.

\subsection*{$\V$-patterns}

If $\CA$ is a small symmetric monoidal $\V$-category, the
$\V$-category of presheaves $\presheaf\CA$ on $\CA$ is a symmetric
monoidal $\V$-category: the convolution of presheaves
$V_1,\dotsc,V_n\in\presheaf{\CA}$ is the $\V$-coend
\begin{equation*}
  V_1 \ast \dotsb \ast V_n = \tint^{A_1,\dotsc,A_n\in\CA}  V_1(A_1)
  \o \dotsb \o V_n(A_n) \o y(A_1\o\dotsb\o A_n) .
\end{equation*}
The Yoneda functor $y\co\CA\to\presheaf\CA$ is a symmetric monoidal
$\V$-functor. (See \citet{Day} and \citet{IK}.)

If $\tau\co\CS\to\CT$ is a symmetric monoidal $\V$-functor between
small symmetric monoidal $\V$-categories, the pull-back functor
$\presheaf{F}$ is a lax symmetric monoidal $\V$-functor; we saw this
in the unenriched case in the introduction, and the proof in the
enriched case is similar.
\begin{definition}
  A \textup{\textbf{$\V$-pattern}} is a symmetric monoidal
  $\V$-functor $\tau:\CS\to\CT$ between small symmetric monoidal
  $\V$-categories $\CS$ and $\CT$ such that
  \begin{equation*}
    \presheaf{\tau}:\presheaf{\CT}\to\presheaf{\CS}
  \end{equation*}
  is a symmetric monoidal $\V$-functor.
\end{definition}

Let $\C$ be a symmetric monoidal $\V$-category which is cocomplete,
and such that the $\V$-functors $A\o B$ preserves colimits in each
variable.  The $\V$-categories of $\tau$-preoperads and $\tau$-operads
in $\C$ are respectively the $\V$-categories
$\PreOp_\tau(\C)=\[\CS,\C\]$ and $\Op_\tau(\C)=\[\CT,\C\]$ of
symmetric monoidal $\V$-functors from $\CS$ and $\CT$ to $\C$.

\subsection*{The monadicity theorem}
We now construct a $\V$-functor $\tau_*$, which generalizes the
functor taking a preoperad to the free operad that it generates.
\begin{proposition}
  \label{adjunction}
  Let $\tau$ be a $\V$-pattern.  The $\V$-adjunction
  \begin{equation*}
    \adjunction{\tau_*}{[\CS,\C]}{[\CT,\C]}{\tau^*}
  \end{equation*}
  induces a $\V$-adjunction between the categories of
  $\tau$-preoperads and $\tau$-operads
  \begin{equation*}
    \adjunction{\tau_*}{\PreOp_\tau(\C)}{\Op_\tau(\C)}{\tau^*} .
  \end{equation*}
\end{proposition}
\begin{proof}
  Let $G$ be a symmetric monoidal $\V$-functor from $\CS$ to $\C$. The
  left Kan $\V$-extension $\tau_*G$ is the $\V$-coend
  \begin{equation*}
    \tau_*G(B) = \tint^{A\in\CS} \CT(\tau A,B) \o G(A) .
  \end{equation*}
  For each $n$, there is a natural $\V$-equivalence
  \begin{multline*}
    \tau_*G(B_1) \o \dotsb \o \tau_*G(B_n) = \textstyle
    \bigotimes_{k=1}^n
    \tint^{A_k\in\CS} \CT(\tau A_k,B_k) \o G( A_k ) \\
    \begin{aligned}
      &\cong \textstyle \tint^{A_1,\dotsc,A_n\in\CS} \bigotimes_k
      \CT(\tau A_k,B_k) \o \bigotimes_{k=1}^n G( A_k ) \\
      & \qquad \qquad \text{since $\o$ preserves $\V$-coends} \\
      &\cong \textstyle \tint^{A_1,\dotsc,A_n\in\CS} \bigotimes_k
      \CT(\tau A_k,B_k) \o G\bigl( \bigotimes_{k=1}^n A_k \bigr) \\
      & \qquad \qquad \text{since $G$ is symmetric monoidal} \\
      &\cong \textstyle \tint^{A_1,\dotsc,A_n\in\CS} \bigotimes_k
      \CT(\tau A_k,B_k) \o \tint^{A\in\CS} \CS\bigl( A ,
      \bigotimes_{k=1}^n A_k \bigr) \o G(A) \\
      & \qquad \qquad \text{by the Yoneda lemma} \\
      &\cong \textstyle \tint^{A\in\CS} \tint^{A_1,\dotsc,A_n\in\CS}
      \bigotimes_{k=1}^n \CT(\tau A_k,B_k) \o \CS\bigl( A ,
      \bigotimes_{k=1}^n A_k \bigr) \o G(A) \\
      & \qquad \qquad \text{by Fubini's theorem for $\V$-coends} \\
      &\cong \textstyle \tint^{A\in\CS} \CT\bigl( \tau A ,
      \bigotimes_{k=1}^n B_k \bigr) \o G(A) \\
      & \qquad \qquad \text{since $\presheaf{\tau}$ is symmetric monoidal} \\
      &= \tau_*G\bigl( \bigotimes_{k=1}^n B_k \bigr) .
    \end{aligned}
  \end{multline*}
  This natural $\V$-equivalence makes $\tau_*G$ into a symmetric
  monoidal $\V$-functor.

  The unit and counit of the $\V$-adjunction between $\tau_*$ and
  $\tau^*$ on $\PreOp_\tau(\C)$ and $\Op_\tau(\C)$ are now induced by
  the unit and counit of the $\V$-adjunction between $\tau_*$ and
  $\tau^*$ on $[\CS,\C]$ and $[\CT,\C]$.
\end{proof}

For the unenriched version of the following result on reflexive
$\V$-coequalizers, see, for example,
\citet[Corollary~1.2.12]{elephant}; the proof in the enriched case is
identical.
\begin{proposition}
  \label{reflexive}
  Let $\C$ be a symmetric monoidal $\V$-category with reflexive
  $\V$-coequalizers. If the tensor product $A\o-$ preserves reflexive
  $\V$-coequalizers, then so does the functor $\SSS(-)$.
\end{proposition}

\begin{corollary}
  \label{create}
  The $\V$-categories $\PreOp_\tau(\C)$ and $\Op_\tau(\C)$ have
  reflexive coequalizers.  
\end{corollary}
\begin{proof}
  The $\V$-coequalizer $\R$ in $[\CS,\C]$ of a reflexive parallel pair
  $\P\reflex^f_g\Q$ in $\[\CS,\C\]$ is computed pointwise: for each
  $X\in\Ob(\CS)$,
  \begin{equation*}
    \P(X) \reflex^{f(X)}_{g(X)} \Q(X) \to \R(X)
  \end{equation*}
  is a reflexive $\V$-coequalizer in $\C$. By
  Proposition~\ref{reflexive}, $\R$ is a symmetric monoidal
  $\V$-functor; thus, $\R$ is the $\V$-coequalizer of the reflexive
  pair $\P\reflex^f_g\Q$ in $\[\CS,\C\]$. The same argument works for
  $\[\CT,\C\]$.
\end{proof}

\begin{proposition}
  If $\tau$ is an essentially surjective $\V$-pattern, the
  $\V$-functor
  \begin{equation*}
    \tau^*:\Op_\tau(\C)\to\PreOp_\tau(\C)
  \end{equation*}
  creates reflexive $\V$-coequalizers.
\end{proposition}
\begin{proof}
  The proof of Corollary \ref{create} shows that the horizontal
  $\V$-functors in the diagram
  \begin{equation*}
    \begin{xy}
      \square/>`-->`>`>/<700,400>[
      \[\CT,\C\]`{[\CT,\C]}`\[\CS,\C\]`{[\CS,\C]};`\tau^*`\tau^*`]
    \end{xy}
  \end{equation*}
  create, and hence preserve, reflexive $\V$-coequalizers. The
  $\V$-functor
  \begin{equation*}
    \tau^*\co[\CT,\C]\to[\CS,\C]
  \end{equation*}
  creates all $\V$-colimits, since $\V$-colimits are computed
  pointwise and $\tau$ is essentially surjective.  It follows that the
  $\V$-functor $\tau^*\co\[\CT,\C\]\to\[\CS,\C\]$ creates reflexive
  $\V$-coequalizers.
\end{proof}

Recall that a $\V$-functor $R\co\CA\to\CB$, with left adjoint
$L\co\CB\to\CA$, is $\V$-monadic if there is an equivalence of
$\V$-categories $\CA\simeq\CB^\TT$, where $\TT$ is the $\V$-monad
associated to the $\V$-adjunction
\begin{equation*}
  \adjunction{L}{\CA}{\CB}{R} .
\end{equation*}
The following is a variant of Theorem II.2.1 of \citet{Dubuc};
reflexive $\V$-coequalizers are substituted for contractible
$\V$-coequalizers, but otherwise, the proof is the same.

\begin{proposition}
  A $\V$-functor $R\co\CA\to\CB$, with left adjoint $L\co\CB\to\CA$,
  is $\V$-monadic if $\CB$ has, and $R$ creates, reflexive
  $\V$-coequalizers.
\end{proposition}

\begin{corollary}
  If $\tau$ is an essentially surjective $\V$-pattern, then the
  $\V$-functor
  \begin{equation*}
    \tau^* \co \Op_\tau(\C) \to \PreOp_\tau(\C)
  \end{equation*}
  is $\V$-monadic.
\end{corollary}

In practice, the $\V$-patterns of interest all have the following
property.
\begin{definition}
  A $\V$-pattern $\tau:\CS\to\CT$ is \textbf{regular} if it is
  essentially surjective and $\CS$ is equivalent to a free symmetric
  monoidal $\V$-category.
\end{definition}

Denote by $\K$ a $\V$-category such that $\CS$ is equivalent to the
free symmetric monoidal $\V$-category $\SSS\K$.  The $\V$-category
$\K$ may be thought of as a generalized set of colours; the theory
associated to a coloured operad is a regular pattern with discrete
$\K$.

\begin{theorem}
  If $\tau$ is a regular pattern and $\C$ is locally finitely
  presentable, then $\Op_\tau(\C)$ is locally finitely presentable.
\end{theorem}
\begin{proof}
  The $\V$-category $\PreOp_\tau(\C)$ is equivalent to the
  $\V$-category $[\K,\C]$, and hence is locally finitely presentable.
  By Lemma~\ref{lift}, $\Op_\tau(\C)$ is cocomplete. The functor
  $\tau_*$ takes finitely presentable objects of $\C$ to finitely
  presentable objects of $\Op_\tau(\C)$, and hence takes finitely
  presentable strong generators of $\PreOp_\tau(\C)$ to finitely
  presentable strong generators of $\Op_\tau(\C)$.
\end{proof}

\section{The simplicial patterns $\Cob$ and $\vCob$} \label{cob}

In this section, we construct simplicial patterns $\Cob=\Cob[d]$,
associated with gluing of oriented $d$-dimensional manifolds along
components of their boundary.

\begin{definition}
  An object $S$ of $\Cob$ is a compact $d$-dimensional manifold, with
  boundary $\p S$ and orientation $o$.
\end{definition}

In the above definition, we permit the manifold $S$ to be
disconnected. In particular, it may be empty. (Note that an empty
manifold of dimension $d$ has a unique orientation.) However, the
above definition should be refined in order to produce a set of
objects of $\Cob$: one way to do this is to add to the data defining
an object of $\Cob$ an embedding into a Euclidean space $\RR^N$,
together with a collared neigbourhood of the boundary $\p S$. We call
this a decorated object.

We now define a simplicial set $\Cob(S,T)$ of morphisms between
objects $S$ and $T$ of $\Cob$.

\begin{definition}
  A hypersurface $\gamma$ in an object $S$ of $\Cob$ consists of a
  closed $(d-1)$-dimensional manifold $M$ together with an embedding
  $\gamma \co M \hookrightarrow S$.
\end{definition}

Given a hypersurface $\gamma$ in $S$, let $S[\gamma]$ be the manifold
obtained by cutting $S$ along the image of $\gamma$. The orientation
$o$ if $S$ induces an orientation $o[\gamma]$ of $S[\gamma]$.

A hypersurface in a decorated object is an open embedding of the
manifold $M\times[-1,1]$ into the complement in $S$ of the collared
neighbourhoods of the boundary. This induces a collaring on the
boundary of $S[\gamma]$. To embed $S[\gamma]$ into $\RR^{N+1}$, we
take the product of the embedding of $S$ into $\RR^N$ and the function
$\chi\circ\pi_2\circ\gamma^{-1}$, where $\pi_2\circ\gamma^{-1}$ is the
function on $S$ equal to the coordinate $t\in[-1,1]$ on the image of
$\gamma$ and undefined elsewhere, and $\chi(t)=\sgn(t)\phi(t)$. Here,
$\phi\in C_c^\infty(-1,1)$ is a non-negative smooth function of
compact support equal to $1$ in a neighbourhood of $0\in(-1,1)$.

A $k$-simplex in the simplicial set $\Cob(S,T)$ of morphisms from $S$
to $T$ consists of the following data:
\begin{enumerate}
\item a closed $(d-1)$-manifold $M$;
\item an isotopy of hypersurfaces, that is a commutative diagram
  \begin{equation*}
    \begin{xy}
      \Vtriangle/@{^{ (}->}`>`>/<500,300>[
      \Delta^k\times M`\Delta^k\times T`\Delta^k;\gamma``]
    \end{xy}
  \end{equation*}
  in which $\gamma$ is an embedding;
\item a fibred diffeomorphism
  \begin{equation*}
    \begin{xy}
      \Vtriangle<500,300>[\Delta^k\times S`{\Delta^k\times T[\gamma]}`\Delta^k;
      \phi``]
    \end{xy}
  \end{equation*}
  compatible with the orientations on its domain and target.
\end{enumerate}
It is straightforward to extend this definition to decorated objects:
the only difference is that $\gamma$ is a fibred open embedding of
$\Delta^k\times M\times[-1,1]$ into $\Delta^k\times T$.

% such that $\phi(t_1,\dots,t_k,x)$ is independent of
% $(t_1,\dots,t_k)\in\Delta^k$ if $x\in\p S$

We now make $\Cob$ into a symmetric monoidal simplicial category. The
composition of $k$-simplices $(M,\gamma,\phi)\in\Cob(S,T)_k$ and
$(N,\delta,\psi)\in\Cob(T,U)_k$ is the $k$-simplex consisting of the
embedding
\begin{equation*}
  \begin{xy}
    \Vtriangle/@{^{ (}->}`>`>/<500,300>[
    \Delta^k\times(M\sqcup N)`\Delta^k\times U`\Delta^k;
    (\psi\circ\gamma)\sqcup\delta``]
  \end{xy}
\end{equation*}
and the fibred diffeomorphism
\begin{equation*}
  \begin{xy}
    \Vtriangle<500,300>[\Delta^k\times S`
    {\Delta^k\times U[(\psi\circ\gamma)\sqcup\delta]}`\Delta^k;\psi\circ\phi``]
  \end{xy}
\end{equation*}
In the special case that the diffeomorphisms $\phi$ and $\psi$ are the
identity, this composition is obtained by taking the union of the
disjoint hypersurfaces $\Delta^k\times M$ and $\Delta^k\times N$ in
$\Delta^k\times U$. It is clear that composition is associative, and
compatible with the face and degeneracy maps between simplices.

The identity $0$-simplex $1_S$ in $\Cob(S,S)$ is associated to the
empty hypersurface in $S$ and the identity diffeomorphism of $S$.

The tensor product of $\Cob$ is simple to describe: it is disjoint
union. When $S_i$ are decorated objects of $\Cob$, $1\le i\le k$,
embedded in $\RR^{N_i}$, we embed $S_1\o\ldots\o S_k$ in
$\RR^{\max(N_i)+1}$ by composing the embedding of $S_i$ with the
inclusion $\RR^{N_i}\hookrightarrow\RR^{\max(N_i)+1}$ defined by
\begin{equation*}
  (t_1,\dots,t_{N_i}) \mapsto (t_1,\dots,t_{N_i},0,\dots,0,i) .
\end{equation*}

Just as modular operads have a directed version, modular dioperads,
the simplicial pattern $\Cob[d]$ has a directed analogue.

\begin{definition}
  An object $S$ of $\vCob$ is a compact $d$-dimensional manifold,
  with orientations $o$ and $\p o$ of $S$ and its boundary $\p S$.
\end{definition}

The boundary $\p S$ of an object $S$ of $\vCob$ is partitioned into
incoming and outgoing parts $\p S=\p_-S\sqcup\p_+S$, according to
whether the boundary is positive or negatively oriented by $\p o$ with
respect to the orientation $o$ of $S$.

\begin{definition}
  A hypersurface $\gamma$ in an object $S$ of $\vCob$ consists of a
  closed $(d-1)$-dimensional oriented manifold $M$ together with an
  embedding $\gamma \co M \hookrightarrow S$.
\end{definition}

Given a hypersurface $\gamma$ in $S$, let $S[\gamma]$ be the manifold
obtained by cutting $S$ along the image of $\gamma$. The orientation
$o$ if $S$ induces an orientation $o[\gamma]$ of $S[\gamma]$, while
the orientations of $\p S$ and $M$ induce an orientation of its
boundary $\p S[\gamma]$.

A $k$-simplex in the simplicial set $\vCob(S,T)$ of morphisms from $S$
to $T$ consists of the following data:
\begin{enumerate}
\item a closed oriented $(d-1)$-manifold $M$;
\item an isotopy of hypersurfaces, that is a commutative diagram
  \begin{equation*}
    \begin{xy}
      \Vtriangle/@{^{ (}->}`>`>/<500,300>[
      \Delta^k\times M`\Delta^k\times T`\Delta^k;\gamma``]
    \end{xy}
  \end{equation*}
  in which $\gamma$ is an embedding;
\item a fibred diffeomorphism
  \begin{equation*}
    \begin{xy}
      \Vtriangle<500,300>[\Delta^k\times S`{\Delta^k\times T[\gamma]}`\Delta^k;
      \phi``]
    \end{xy}
  \end{equation*}
  compatible with the orientations on its domain and target.
\end{enumerate}

We may make $\vCob$ into a symmetric monoidal simplicial category in
the same way as for $\Cob$. The definition of decorated objects and
decorated morphisms is also easily extended to this setting.

Let $\cob$ be the simplicial groupoid of connected oriented surfaces
and their oriented diffeomorphisms, with skeleton
\begin{equation*}
  \bigsqcup_{g,n} \, \Diff_+(S_{g,n}) .
\end{equation*}
The embedding $\cob\hookrightarrow\Cob$ extends to an essentially
surjective $\V$-functor $\SSS\cob\to\Cob$. Similarly, if $\vcob$ is
the simplicial groupoid of connected oriented surfaces with oriented
boundary and their oriented diffeomorphisms, the $\V$-functor
$\SSS\vcob\to\vCob$ is essentially surjective.

We can now state the main theorem of this section.
\begin{theorem}
  The simplicial functors $\cob\hookrightarrow\Cob$ and
  $\vcob\hookrightarrow\vCob$ induce regular simplicial patterns
  $\SSS\cob\to\Cob$ and $\SSS\vcob\to\vCob$
\end{theorem}

A $d$-dimensional modular operad, respectively dioperad, is an operad
for the simplicial pattern $\SSS\cob\to\Cob$, respectively
$\SSS\vcob\to\vCob$.

\subsection*{One-dimensional modular operads}

When $d=1$, the category $\cob$ has a skeleton with two objects, the
interval $I$ and the circle $S$. The simplicial group
$\cob(I,I)\cong\Diff_+[0,1]$ is contractible and the simplicial group
$\cob(S,S)$ is homotopy equivalent to $\SO(2)$.

The simplicial set $\Cob(I^{\o k},I)$ is empty for $k=0$ and
contractible for each $k>0$, and the simplicial set $\Cob(I^{\o k},S)$
is empty for $k=0$ and homotopy equivalent to $\SO(2)$ for each $k>0$.
Thus, a 1-dimensional modular operad in a symmetric monoidal category
$\C$ consists of a homotopy associative algebra $A$, and a homotopy
trace from $A$ to a homotopy $\SO(2)$-module $M$. If $\C$ is discrete,
a 1-dimensional modular operad is simply a non-unital associative
algebra $A$ in $\C$ together with an object $M$ in $\C$ and a trace
$\tr\co A\to M$.

\subsection*{One-dimensional modular dioperads}

When $d=1$, the category $\vcob$ has a skeleton with five objects, the
intervals $I_-^-$, $I_-^+$, $I_+^-$ and $I_+^+$, representing the
$1$-manifold $[0,1]$ with the four different orientations of its
boundary, and the circle $S$.  The simplicial groups $\vcob(I_a,I_a)$
are contractible, and the simplicial group $\vcob(S,S)$ is homotopy
equivalent to $\SO(2)$.

A 1-dimensional modular dioperad $\P$ in a discrete symmetric monoidal
category $\C$ consists of associative algebras $A=\P(I_-^+)$ and
$B=\P(I_+^-)$ in $\C$, a $(A,B)$-bimodule $Q=\P(I_-^-)$, a
$(B,A)$-bimodule $R=\P(I_+^+)$, and an object $M=\P(S)$, together with
morphisms
\begin{align*}
  \alpha & \co Q \o_B R \to A , & \beta & \co R \o_A Q \to B , &
  \tr_A & \co A \to M , & \tr_B & \co B \to M .
\end{align*}
Denote the left and right actions of $A$ and $B$ on $Q$ and $R$ by
$\lambda_Q\co A\o Q\to Q$, $\rho_Q\co Q\o B\to Q$, $\lambda_R\co B\o R\to R$
and $\rho_R\co R\o A\to R$. The above data must in addition satisfy
the following conditions:
\begin{itemize}
\item $\alpha$ and $\beta$ are morphisms of $(A,A)$-bimodules and
  $(B,B)$-bimodules respectively;
\item $\tr_A$ and $\tr_B$ are traces;
\item $\lambda_Q\circ(\alpha\o Q)=\rho_Q\circ(R\o\beta)\co Q\o_BR\o_A
  Q\to Q$ and $\lambda_R\circ(\beta\o R) =
  \rho_R\circ(Q\o\alpha)\co R\o_AQ\o_B R\to R$;
\item $\tr_A\circ\alpha\co Q\o_BR\to M$ and $\tr_B\circ\beta\co R\o_AQ\to M$
  are equal on the isomorphic objects $(Q\o_BR)\o_{A^\op\o A}\1$ and
  $(R\o_AQ)\o_{B^\op\o B}\1$.
\end{itemize}
That is, a 1-dimensional modular dioperad is the same thing as a
\textbf{Morita context} (\citet{Morita}) $(A,B,P,Q,\alpha,\beta)$ with
trace $(\tr_A,\tr_B)$.

\subsection*{Two-dimensional topological field theories}

Let $\CG$ be the discrete pattern introduced in Section~1 whose operads
are modular operads. There is a natural morphism of patterns
\begin{equation*}
  \alpha\co\Cob[2]\to\CG ;
\end{equation*}
thus, application of $\alpha^*$ to a modular operad gives rise to a
2-dimensional modular operad. But modular operads and what we call
2-dimensional modular operads are quite different: a modular preoperad
$\CP$ consists of a sequence of $\SS_n$-modules $\CP\(g,n\)$, while a
2-dimensional modular preoperad is a $\cob$-module.

An example of a $2$-dimensional modular operad is the terminal one,
for which $\CP(S)$ is the unit $\1$ for each surface $S$. A more
interesting one comes from conformal field theory: the underlying
2-dimensional modular operad associates to an oriented surface with
boundary $S$ the moduli space $\CN(S)$ of conformal structures on $S$.
(More accurately, $\CN(S)$ is the simplicial set whose $n$-simplices
are the $n$-parameter smooth families of conformal structures
parametrized by the $n$-simplex.) The space $\CN(S)$ is contractible
for all $S$, since it is the space of smooth sections of a fibre
bundle over $M$ with contractible fibres. To define the structure of a
2-dimensional modular operad on $\CN$ amounts to showing that
conformal structures may be glued along circles: this is done by
choosing a Riemannian metric in the conformal class which is flat in a
neighbourhood of the boundary, such that the boundary is geodesic, and
each of its components has length $1$.

A bilinear form $t\co M\o M \to\1$ on a cochain complex where $M$ is
\textbf{non-degenerate} if it induces a quasi-isomorphism between $M$
and $\Hom(M,\1)$.

\begin{definition}
  Let $\P$ be a 2-dimensional modular operad.  A $\P$-algebra is a
  morphism of 2-dimensional modular operads
  \begin{equation*}
    \rho \co \P \to \alpha^*\End(M,t) ,
  \end{equation*}
  where $M$ is a cochain complex with non-degenerate bilinear form
  $t\co M\o M \to\1$.

  A topological conformal field theory is a $C_*(\CN)$-algebra.
\end{definition}

In the theory of infinite loop spaces, one defines an
$E_\infty$-algebra as an algebra for an operad $\E$ such that $\E(n)$
is contractible for all $n$. Similarly, as shown by \citet{braided},
an $E_2$-algebra is an algebra for a braided operad $\E$ such that
$\E(n)$ is contractible for all $n$. Motivated by this, we make the
following definition.
\begin{definition}
  A $2$-dimensional topological field theory is a pair consisting of a
  $2$-dimensional modular operad $\E$ in the category of cochain
  complexes such that $\E(S)$ is quasi-isomorphic to $\1$ for all
  surfaces $S$, and an $\E$-algebra $(M,t)$.
\end{definition}
In particular, a topological conformal field theory is a
$2$-dimensional topological field theory.

\section*{Appendix. Enriched categories}

In this appendix, we recall some results of enriched category theory.
Let $\V$ be a closed symmetric monoidal category; that is, $\V$ is a
symmetric monoidal category such that the tensor product functor $-\o
Y$ from $\V$ to itself has a right adjoint for all objects $Y$,
denoted $[Y,-]$: in other words,
\begin{equation*}
  [X\o Y,Z] \cong [X,[Y,Z]] .
\end{equation*}
Throughout this paper, we assume that $\V$ is complete and cocomplete.
Denote by $A\mapsto A_0$ the continuous functor $A_0=\V(\1,A)$ from
$\V$ to $\Set$, where $\1$ is the unit of $\V$.

Let $\VCat$ be the 2-category whose objects are $\V$-categories, whose
1-mor\-phisms are $\V$-functors, and whose 2-morphisms are
$\V$-natural transformations. Since $\V$ is closed, $\V$ is itself a
$\V$-category.

Applying the functor $(-)_0\co\V\to\Set$ to a $\V$-category $\C$, we
obtain its underlying category $\C_0$; in this way, we obtain a
2-functor
\begin{equation*}
  (-)_0 \co \VCat \to \Cat .
\end{equation*}

Given a small $\V$-category $\CA$, and a $\V$-category $\CB$, there is
a $\V$-category $[\CA,\CB]$, whose objects are the $\V$-functors from
$\CA$ to $\CB$, and such that
\begin{equation*}
  [\CA,\CB](F,G) = \tint_{A\in\CA} \CB(FA,GA) .  
\end{equation*}
There is an equivalence of categories $\VCat(\CA,\CB)\simeq[\CA,\CB]_0$.

\subsection*{The Yoneda embedding for $\V$-categories}
The opposite of a $\V$-category $\C$ is the $\V$-category $\C^\op$
with
\begin{equation*}
  \C^\op(A,B) = \C(B,A) .
\end{equation*}
If $\CA$ is a small $\V$-category, denote by $\presheaf\CA$ the
$\V$-category of presheaves
\begin{equation*}
  \presheaf{\CA} = [\CA^\op,\V] .
\end{equation*}
By the $\V$-Yoneda lemma, there is a full, faithful $\V$-functor
$y\co\CA\to\presheaf{\CA}$, with
\begin{equation*}
  y(A) = \CA(-,A) , \quad A\in\Ob(\CA) .
\end{equation*}
If $F\co\CA\to\CB$ is a $\V$-functor, denote by
$\presheaf{F}\co\presheaf{\CB}\to\presheaf{\CA}$ the $\V$-functor
induced by $F$.

\subsection*{Cocomplete $\V$-categories}
A $\V$-category $\C$ is \textbf{tensored} if there is a $\V$-functor
$\o\co\V\o\C\to\C$ together with a $\V$-natural equivalence of
functors $\C(X\o A,B) \cong [X,\C(A,B)]$ from $\V^\op\o\C^\op\o\C$ to
$\V$. For example, the $\V$-category $\V$ is itself tensored.

Let $\CA$ be a small $\V$-category, let $F$ be a $\V$-functor from
$\CA^\op$ to $\V$, and let $G$ be a $\V$-functor from $\CA$ to a
$\V$-category $\C$. If $B$ is an object of $\C$, denote by
$\<G,B\>:\C^\op\to\V$ the presheaf such that $\<G,C\>(A)=[GA,C]$. The
\textbf{weighted colimit} of $G$, with weight $F$, is an object
$\tint^{A\in\CA} FA \o GA$ of $\C$ such that there is a natural
isomorphism
\begin{equation*}
 \presheaf{\C}( F , \< G , - \> ) \cong \C \bigl( \tint^{A\in\CA} FA
 \o GA , - \bigr) .
\end{equation*}
If $\C$ is tensored, $\tint^{A\in\CA} F(A)\o G(A)$ is the coequalizer
of the diagram
\begin{equation*}
  \bigsqcup_{A_0,A_1\in\Ob(\CA)} \CA(A_0,A_1) \o FA_1 \o GA_0
  \two \bigsqcup_{A\in\Ob(\CA)} FA \o GA ,
\end{equation*}
where the two morphisms are induced by the action on $F$ and coaction
on $G$ respectively.

A $\V$-category $\C$ is \textbf{cocomplete} if it has all weighted
colimits, or equivalently, if it satisfies the following conditions:
\begin{enumerate}
\item the category $\C_0$ is cocomplete;
\item for each object $A$ of $\C$, the functor $\C(-,A) \co \C_0 \to \V$
  transforms colimits into limits;
\item $\C$ is tensored.
\end{enumerate}
In particular, $\V$-categories $\presheaf{\CA}$ of presheaves are
cocomplete.

Let $\C$ be a cocomplete $\V$-category. A $\V$-functor $F\co\CA\to\CB$
between small $\V$-categories gives rise to a $\V$-adjunction
\begin{equation*}
  \adjunction{F_*}{[\CA,\C]}{[\CB,\C]}{F^*} ,
\end{equation*}
that is, an adjunction in the 2-category $\VCat$. The functor $F_*$
is called the (pointwise) left $\V$-Kan extension of along $F$; it is
the $\V$-coend
\begin{equation*}
  F_*G(-) = \tint^{A\in\CA} \CA(FA,-) \o GA .
\end{equation*}

\subsection*{Cocomplete categories of $\V$-algebras}

A $\V$-monad $\TT$ on a $\V$-category $\C$ is a monad in the full
sub-2-category of $\VCat$ with unique object $\C$. If $\C$ is small,
this is the same thing as a monoid in the monoidal category
$\VCat(\C,\C)$.

The following is the enriched version of a result of \citet{Linton},
and is proved in exactly the same way.
\begin{lemma}
  \label{lift}
  Let $\TT$ be a $\V$-monad on a cocomplete $\V$-category $\C$ such
  that the $\V$-category of algebras $\C^\TT$ has reflexive
  $\V$-coequalizers. Then the $\V$-category of $\TT$-algebras $\C^\TT$
  is cocomplete.
\end{lemma}
\begin{proof}
  We must show that $\C^\TT$ has all weighted colimits. Let $\CA$ be a
  small $\V$-category, let $F$ be a weight, and let $G:\CA\to\C^\TT$
  be a diagram of $\TT$-algebras. Then the weighted colimit
  $\tint^{A\in\CA} F(A) \o G(A)$ is a reflexive coequalizer
  \begin{equation*}
    \TT \bigl( \tint^{A\in\CA} FA \o \TT RGA \bigr) \reflex \TT
    \bigl( \tint^{A\in\CA} FA \o RGA \bigr) .
  \qedhere
  \end{equation*}
\end{proof}

\subsection*{Locally finitely presentable $\V$-categories}

In studying algebraic theories using categories, locally finite
presentable categories plays a basic role. When the closed symmetric
monoidal categery $\V$ is locally finitely presentable, with finitely
presentable unit $\1$, these have a generalization to enriched category
theory over $\V$, due to \citet{enriched}.  (There is a more general
theory of locally presentable categories, where the cardinal
$\aleph_0$ is replaced by an arbitrary regular cardinal; this
extension is straightforward, but we do not present it here in order
to simplify exposition.)

Examples of locally finitely presentable closed symmetric monoidal
categories include the categories of sets, groupoids, categories,
simplicial sets, and abelian groups --- the finitely presentable
objects are respectively finite sets, groupoids and categories,
simplicial sets with a finite number of nondegenerate simplices, and
finitely presentable abelian groups. A less obvious example is the
category of symmetric spectra of \citet{HSS}.

An object $A$ in a $\V$-category $\C$ is \textbf{finitely presentable}
if the functor
\begin{equation*}
  \C(A,-) : \C \to \V
\end{equation*}
preserves filtered colimits. A \textbf{strong generator} in a
$\V$-category $\C$ is a set $\{G_i\in\Ob(\C)\}_{i\in I}$ such that the
functor
  \begin{equation*}
    A \mapsto \bigotimes_{i\in I} \C(G_i,A) : \C \mapsto \V^{\o I}
  \end{equation*}
  reflects isomorphisms.
\begin{definition}
  A \textbf{locally finitely presentable} $\V$-category $\C$ is a
  cocomplete $\V$-category with a finitely presentable strong
  generator (i.e.\ the objects making up the strong generator are
  finitely presentable).
\end{definition}

\subsection*{Pseudo algebras}

A 2-monad $\TT$ on a 2-category $\CC$ is by definition a $\Cat$-monad
on $\CC$.  Denote the composition of the 2-monad $\TT$ by
$m\co\TT\TT\to\TT$, and the unit by $\eta\co1\to\TT$.

Associated to a 2-monad $\TT$ is the 2-category $\CC^\TT$ of pseudo
$\TT$-algebras (\citet{lax} and \citet{pseudo}; see also
\citet{Marmolejo}). A pseudo $\TT$-algebra is an object $\CA$ of
$\CC$, together with a morphism $a\co\TT\CA\to\CA$, the composition,
and invertible 2-morphisms
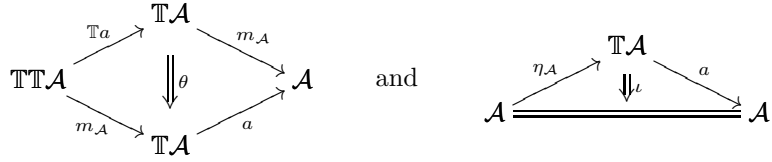
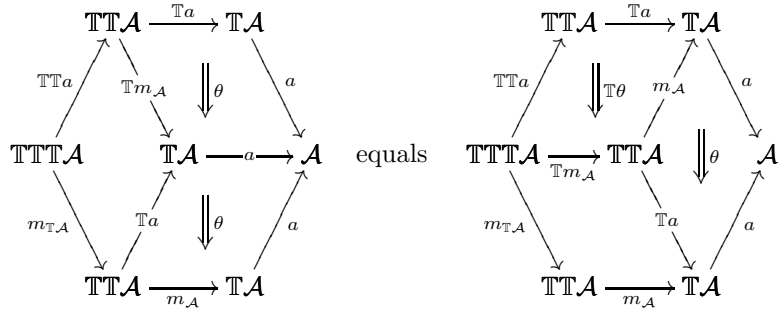
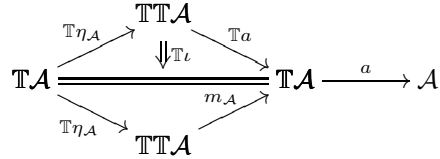
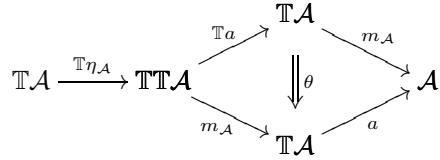
\begin{align*}
  \begin{aligned}
    \begin{xy}
      \morphism(0,0)<500,250>[\TT\TT\CA`\TT\CA;\TT a]
      \morphism(0,0)|b|<500,-250>[\TT\TT\CA`\TT\CA;m_\CA]
      \morphism(500,250)<500,-250>[\TT\CA`\CA;m_\CA]
      \morphism(500,-250)|b|<500,250>[\TT\CA`\CA;a]
      \morphism(500,100)|r|/=>/<0,-200>[`;\theta]
    \end{xy}
  \end{aligned}
  & & \text{and} & &
  \begin{aligned}
    \begin{xy}
      \morphism(0,0)<500,250>[\CA`\TT\CA;\eta_\CA]
      \morphism(0,0)/=/<1000,0>[\CA`\CA;]
      \morphism(500,250)<500,-250>[\TT\CA`\CA;a]
      \morphism(500,150)|r|/=>/<0,-100>[`;\iota]
    \end{xy}
  \end{aligned}
\end{align*}
such that
\begin{gather*}
  \begin{aligned}
    \begin{xy}
      \morphism(0,0)<250,500>[\TT\TT\TT\CA`\TT\TT\CA;\TT\TT a]
      \morphism(0,0)|l|<250,-500>[\TT\TT\TT\CA`\TT\TT\CA;m_{\TT\CA}]
      \morphism(250,500)|m|<250,-500>[\TT\TT\CA`\TT\CA;\TT m_\CA]
      \morphism(250,-500)|m|<250,500>[\TT\TT\CA`\TT\CA;\TT a]
      \morphism(250,500)<500,0>[\TT\TT\CA`\TT\CA;\TT a]
      \morphism(500,0)|m|<500,0>[\TT\CA`\CA;a]
      \morphism(250,-500)|b|<500,0>[\TT\TT\CA`\TT\CA;m_\CA]
      \morphism(750,500)<250,-500>[\TT\CA`\CA;a]
      \morphism(750,-500)|b|<250,500>[\TT\CA`\CA;a]
      \morphism(600,350)|r|/=>/<0,-200>[`;\theta]
      \morphism(600,-150)|r|/=>/<0,-200>[`;\theta]
    \end{xy}
  \end{aligned}
  \quad \text{equals} \quad
  \begin{aligned}
    \begin{xy}
      \morphism(0,0)<250,500>[\TT\TT\TT\CA`\TT\TT\CA;\TT\TT a]
      \morphism(0,0)|l|<250,-500>[\TT\TT\TT\CA`\TT\TT\CA;m_{\TT\CA}]
      \morphism(0,0)|b|<500,0>[\TT\TT\TT\CA`\TT\TT\CA;\TT m_\CA]
      \morphism(250,500)<500,0>[\TT\TT\CA`\TT\CA;\TT a]
      \morphism(500,0)|m|<250,500>[\TT\TT\CA`\TT\CA;m_\CA]
      \morphism(500,0)|m|<250,-500>[\TT\TT\CA`\TT\CA;\TT a]
      \morphism(250,-500)|b|<500,0>[\TT\TT\CA`\TT\CA;m_\CA]
      \morphism(750,500)<250,-500>[\TT\CA`\CA;a]
      \morphism(750,-500)|b|<250,500>[\TT\CA`\CA;a]
      \morphism(350,350)|r|/=>/<0,-200>[`;\TT\theta]
      \morphism(750,100)|r|/=>/<0,-200>[`;\theta]
    \end{xy}
  \end{aligned}
\end{gather*}
and
\begin{gather*}
  \begin{aligned}
    \begin{xy}
      \morphism(0,0)<500,250>[\TT\CA`\TT\TT\CA;\TT\eta_\CA]
      \morphism(0,0)|b|<500,-250>[\TT\CA`\TT\TT\CA;\TT\eta_\CA]
      \morphism(0,0)/=/<1000,0>[\TT\CA`\TT\CA;]
      \morphism(500,250)<500,-250>[\TT\TT\CA`\TT\CA;\TT a]
      \morphism(500,-250)<500,250>[\TT\TT\CA`\TT\CA;m_\CA]
      \morphism(1000,0)<500,0>[\TT\CA`\CA;a]
      \morphism(500,150)|r|/=>/<0,-100>[`;\TT\iota]
    \end{xy}
  \end{aligned}
  \intertext{equals}
  \begin{aligned}
    \begin{xy}
      \morphism(-500,0)<500,0>[\TT\CA`\TT\TT\CA;\TT\eta_\CA]
      \morphism(0,0)<500,250>[\TT\TT\CA`\TT\CA;\TT a]
      \morphism(0,0)|b|<500,-250>[\TT\TT\CA`\TT\CA;m_\CA]
      \morphism(500,250)<500,-250>[\TT\CA`\CA;m_\CA]
      \morphism(500,-250)|b|<500,250>[\TT\CA`\CA;a]
      \morphism(500,100)|r|/=>/<0,-200>[`;\theta]
    \end{xy}
  \end{aligned}
\end{gather*}

The lax morphisms of $\CC^\TT$ are pairs $(f,\phi)$ consisting of a
morphism $f\co\CA\to\CB$ and a 2-morphism
\begin{equation*}
  \begin{xy}
    \morphism(0,0)<500,250>[\TT\CA`\TT\CB;\TT f]
    \morphism(0,0)|b|<500,-250>[\TT\CA`\CA;a]
    \morphism(500,250)<500,-250>[\TT\CB`\CB;b]
    \morphism(500,-250)|b|<500,250>[\CA`\CB;f]
    \morphism(500,100)|r|/=>/<0,-200>[`;\phi]
  \end{xy}
\end{equation*}
such that
\begin{gather*}
  \begin{aligned}
    \begin{xy}
      \morphism(0,0)<250,500>[\TT\TT\CA`\TT\TT\CB;\TT\TT f]
      \morphism(0,0)|l|<250,-500>[\TT\TT\CA`\TT\CA;m_\CA]
      \morphism(250,500)|m|<250,-500>[\TT\TT\CB`\TT\CB;m_\CB]
      \morphism(250,-500)|m|<250,500>[\TT\CA`\TT\CB;\TT f]
      \morphism(250,500)<500,0>[\TT\TT\CB`\TT\CB;\TT b]
      \morphism(500,0)|m|<500,0>[\TT\CB`\CB;b]
      \morphism(250,-500)|b|<500,0>[\TT\CA`\CA;a]
      \morphism(750,500)<250,-500>[\TT\CB`\CB;b]
      \morphism(750,-500)|b|<250,500>[\CA`\CB;f]
      \morphism(600,350)|r|/=>/<0,-200>[`;\theta]
      \morphism(600,-150)|r|/=>/<0,-200>[`;\phi]
    \end{xy}
  \end{aligned}
  \quad \text{equals} \quad
  \begin{aligned}
    \begin{xy}
      \morphism(0,0)<250,500>[\TT\TT\CA`\TT\TT\CB;\TT f]
      \morphism(0,0)|l|<250,-500>[\TT\TT\CA`\TT\CA;m_\CA]
      \morphism(0,0)|m|<500,0>[\TT\TT\CA`\TT\CA;\TT a]
      \morphism(250,500)<500,0>[\TT\TT\CB`\TT\CB;\TT b]
      \morphism(500,0)|m|<250,500>[\TT\CA`\TT\CB;\TT f]
      \morphism(500,0)|m|<250,-500>[\TT\CA`\CA;a]
      \morphism(250,-500)|b|<500,0>[\TT\CA`\CA;a]
      \morphism(750,500)<250,-500>[\TT\CB`\CB;b]
      \morphism(750,-500)|b|<250,500>[\CA`\CB;f]
      \morphism(350,350)|r|/=>/<0,-200>[`;\TT\phi]
      \morphism(350,-150)|r|/=>/<0,-200>[`;\theta]
      \morphism(750,100)|r|/=>/<0,-200>[`;\theta]
    \end{xy}
  \end{aligned}
\end{gather*}
and
\begin{gather*}
  \begin{aligned}
    \begin{xy}
      \morphism(-500,250)<500,250>[\CA`\TT\CA;\eta_\CA]
      \morphism(0,500)<500,-250>[\TT\CA`\TT\CB;\TT f]
      \morphism(-500,250)|b|<500,-250>[\CA`\CB;f]
      \morphism(0,0)|m|<500,250>[\CB`\TT\CB;\eta_\CB]
      \morphism(0,0)/=/<1000,0>[\CB`\CB;]
      \morphism(500,250)<500,-250>[\TT\CB`\CB;b]
      \morphism(500,150)|r|/=>/<0,-100>[`;\iota]
    \end{xy}
  \end{aligned}
  \intertext{equals}
  \begin{aligned}
    \begin{xy}
      \morphism(-500,-250)<500,250>[\CA`\TT\CA;\eta_\CA]
      \morphism(-500,-250)/=/<1000,0>[\CA`\CA;]
      \morphism(0,0)<500,250>[\TT\CA`\TT\CB;\TT f]
      \morphism(0,0)|m|<500,-250>[\TT\CA`\CA;a]
      \morphism(500,250)<500,-250>[\TT\CB`\CB;b]
      \morphism(500,-250)|b|<500,250>[\CA`\CB;f]
      \morphism(500,100)|r|/=>/<0,-200>[`;\phi]
      \morphism(0,-100)|r|/=>/<0,-100>[`;\iota]
    \end{xy}
  \end{aligned}
\end{gather*}
A lax morphism $(f,\phi)$ is a morphism if the 2-morphism $\phi$ is
invertible.

The 2-morphisms $\gamma\co(f,\phi)\to(\tilde{f},\tilde{\phi})$ of
$\CC^\TT$ are 2-morphisms $\gamma\co f\to\tilde{f}$ such that
\begin{align} \label{2morphism}
  \begin{aligned}
    \begin{xy}
      \morphism(0,0)<500,250>[\TT\CA`\TT\CB;\TT f]
      \morphism(0,0)|b|<500,-250>[\TT\CA`\CA;a]
      \morphism(500,250)<500,-250>[\TT\CB`\CB;b]
      \morphism(500,-250)|a|/{@{>}@/^1em/}/<500,250>[\CA`\CB;f]
      \morphism(500,-250)|b|/{@{>}@/_1em/}/<500,250>[\CA`\CB;\tilde{f}]
      \morphism(500,100)|r|/=>/<0,-200>[`;\phi]
      \morphism(750,-75)|r|/=>/<0,-100>[`;\gamma]
    \end{xy}
  \end{aligned}
  & & \text{equals} & &
  \begin{aligned}
    \begin{xy}
      \morphism(0,0)|a|/{@{>}@/^1em/}/<500,250>[\TT\CA`\TT\CB;\TT f]
      \morphism(0,0)|b|/{@{>}@/_1em/}/<500,250>[\TT\CA`\TT\CB;\TT\tilde{f}]
      \morphism(0,0)|b|<500,-250>[\TT\CA`\CA;a]
      \morphism(500,250)<500,-250>[\TT\CB`\CB;a]
      \morphism(500,-250)|b|<500,250>[\CA`\CB;\tilde{f}]
      \morphism(500,100)|r|/=>/<0,-200>[`;\tilde{\phi}]
      \morphism(250,175)|r|/=>/<0,-100>[`;\TT\gamma]
    \end{xy}
  \end{aligned}
\end{align}

\subsection*{$\V$-categories of pseudo algebras}

If $\CA$ and $\CB$ are pseudo $\TT$-algebras, and the underlying
category of $\CA$ is small, we saw that there is a $\V$-category
$[\CA,\CB]$ whose objects are $\V$-functors $f\co\CA\to\CB$ between
the underlying $\V$-categories.  Using $[\CA,\CB]$, we now define a
$\V$-category $\[\CA,\CB\]$ whose objects are morphisms
$f\co\CA\to\CB$ of pseudo $\TT$-algebras.  If $f_0$ and $f_1$ are
morphisms of pseudo $\TT$-algebras, $\[\CA,\CB\](f_0,f_1)$ is defined
as the equalizer
\begin{equation*}
  \[\CA,\CB\](f_0,f_1) \too^{\psi_{f_0f_1}} [\CA,\CB](f_0,f_1)
  \two^{\phi_1(b\circ\TT-)}_{(-\circ a)\phi_0}
  [\TT\CA,\CB](b \circ\TT f_0 , f_1\circ a)
\end{equation*}
Of course, this is the internal version of \eqref{2morphism}: the
2-morphisms $\gamma\co f_0\to f_1$ of $\CC^\TT$ are the elements of the
set $|\[\CA,\CB\](f_0,f_1)|$.

The composition morphism
\begin{equation*}
  \[\CA,\CB\](f_0,f_1)\o\[\CA,\CB\](f_1,f_2)
  \too^{m^{\[\CA,\CB\]}_{f_0f_1f_2}} \[\CA,\CB\](f_0,f_2)
\end{equation*}
is the universal arrow for the coequalizer $\[\CA,\CB\](f_0,f_2)$,
whose existence is guaranteed by the commutativity of the diagram
\begin{multline*}
  \[\CA,\CB\](f_0,f_1)\o\[\CA,\CB\](f_1,f_2)
  \too^{\psi_{f_0f_1}\o\psi_{f_1f_2}}
  [\CA,\CB](f_0,f_1)\o[\CA,\CB](f_1,f_2) \\
  \too^{m^{[\CA,\CB]}_{f_0f_1f_1}} [\CA,\CB](f_0,f_2)
  \two^{\phi_2(b\circ\TT-)}_{(-\circ a)\phi_0}
  [\TT\CA,\CB](b\circ\TT f_0 , f_2\circ c)
\end{multline*}
Indeed, we have
\begin{multline*}
  \phi_2(b\circ\TT-) \* m^{[\CA,\CB]}_{f_0f_1f_1} \*
  (\psi_{f_0f_1}\o\psi_{f_1f_2}) \\
  \begin{aligned}
    &= m^{[\TT\CA,\CB]}_{b\circ\TT f_0,b\circ\TT f_1,f_2\circ a} \* (
    (b\circ\TT-)\*\psi_{f_0f_1} \o
    \phi_2(b\circ\TT-)\*\psi_{f_1f_2} ) \\
    &= m^{[\TT\CA,\CB]}_{b\circ\TT f_0,b\circ\TT f_1,f_2\circ a} \* (
    (b\circ\TT-)\*\psi_{f_0f_1} \o
    (-\circ a)\phi_1\*\psi_{f_1f_2} ) \\
    &= m^{[\TT\CA,\CB]}_{b\circ\TT f_0,f_1\circ a,f_2\circ a} \* (
    \phi_1(b\circ\TT-)\*\psi_{f_0f_1} \o (-\circ a)\*\psi_{f_1f_2} ) \\
    &= m^{[\TT\CA,\CB]}_{b\circ\TT f_0,f_1\circ a,f_2\circ a} \* (
    (-\circ a)\phi_0\*\psi_{f_0f_1} \o (-\circ a)\*\psi_{f_1f_2} ) \\
    &= (-\circ a)\phi_0 \* m^{[\CA,\CB]}_{f_0f_1f_1} \*
    (\psi_{f_0f_1}\o\psi_{f_1f_2}) .
  \end{aligned}
\end{multline*}
Associativity of this composition is proved by a calculation along the
same lines for the iterated composition map
\begin{equation*}
  \[\CA,\CB\](f_0,f_1)\o\[\CA,\CB\](f_1,f_2)\o\[\CA,\CB\](f_2,f_3)
  \to \[\CA,\CB\](f_0,f_3) .
\end{equation*}
The unit $v_f\co\1\to\[\CA,\CB\](f,f)$ of the $\V$-category
$\[\CA,\CB\]$ is the universal arrow for the coequalizer
$\[\CA,\CB\](f,f)$, whose existence is guaranteed by the
commutativity of the diagram
\begin{equation*}
  \1 \too^{u_f} [\CA,\CB](f,f)
  \two^{\phi(b\circ\TT-)}_{(-\circ a)\phi}
  [\TT\CA,\CB](b\circ\TT f , f\circ a)
\end{equation*}
Indeed, both $\phi(b\circ\TT-)\*u_f$ and $(-\circ a)\phi\*u_f$ are
equal to
\begin{equation*}
  \phi\in\V(\1,[\TT\CA,\CB](b\circ\TT f , f\circ a)) .
\end{equation*}

In a remark at the end of Section~3 of \citep{IK}, Im and Kelly
dismiss the existence of the $\V$-categories of pseudo algebras. The
construction presented here appears to get around their objections.

\bibliography{moduli}

\bibliographystyle{plainnat}

\end{document}